\documentclass[12pt]{article}
\title{ An Introduction to Finite Fibonomial Calculus}
\author{Ewa Krot \\
\\Institute of Computer Science, Bia{\l}ystok University\\
PL-15-887 Bia{\l}ystok, ul.Sosnowa 64, POLAND\\
e-mail: ewakrot@wp.pl}

\usepackage{amsmath,amsthm}

\chardef\bslash=`\\ 
\newtheorem{ex}{Example}[section]
\newtheorem{defn}{Definition}[section]

\newtheorem{rem}{Remark}[section]
\newtheorem{thm}{Theorem}[section]
\newtheorem{prop}{Proposition}[section]
\newtheorem{cor}{Corollary}[section]

\begin{document}
\maketitle
\begin{abstract}
This is an indicatory presentation of main definitions and theorems of Fibonomial
 Calculus which is a special case of $\psi$-extented Rota's finite
 operator calculus \cite{1}.
\end{abstract}

\section{Fibonomial coefficients}
The famous Fibonacci sequence$\{F_{n}\}_{n \geq 0}$
$$\Bigg\{\begin{array}{l}F_{n+2}=F_{n+1}+F_{n}\\
F_{0}=0,\;F_{1}=1\end{array}$$ is attributed and refered to the
first edition (lost) of "Liber Abaci" (1202) by Leonardo Fibonacci
(Pisano)(see edition from 1228 reproduced as "Il Liber Abaci di
Leonardo Pisano publicato secondo la lezione Codice Maglibeciano
by Baldassarre Boncompagni in Scritti di Leonardo Pisano" , vol.
1,(1857)Rome).

In order to specify what a "Fibonomial Calculus" is let us define
for thr sequence $F=\left\{F_{n}\right\}_{n \geq 0}$ what
 follows:
\renewcommand{\labelenumi}{(\arabic{enumi})}
\begin{enumerate}
\item $F$-factorial:
$$F_{n}!=F_{n}F_{n-1}...F_{2}F_{1},\;\;\;F_{0}!=1.$$
\item $F$-binomial (Fibonomial ) coefficients \cite{4}:
$$\binom{n}{k}_{F}=\frac{n_{F}^{\underline{k}}}{k_{F}!}=
\frac{F_{n}F_{n-1}\ldots F_{n-k+1}}{F_{k}F_{k-1}\ldots F_{2}F_{1}}=
\frac{F_{n}!}{F_{k}!F_{n-k}!},\;\;\;\binom{n}{0}_{F}=1.$$
\end{enumerate}
Some properties of $\binom{n}{k}_{F}$ are:
\renewcommand{\labelenumi}{(\alph{enumi})}
\begin{enumerate}
\item $\binom{n}{k}_{F}=\binom{n}{n-k}_{F}$, (symmetry);
\item $F_{n-k}\binom{n}{k}_{F}=F_{n}\binom{n-1}{k}_{F}$;
\item $\binom{n}{k}_{F}\in {\bf N}$ for every $n,k \in {\bf N}\cup {0}$.
\end{enumerate}
\section{Operators and polynomial sequences}
Let {\bf P} be the algebra of polynomials over the field {\bf K}
of characteristic zero.
\begin{defn}
 The linear operator   $\partial_{F}:{\bf P}\rightarrow{\bf P}$   such that
   $\partial_{F}x^{n}=F_{n}x^{n-1}$  for $n \geq 0$ is named the
   $F$-derivative.
\end{defn}
\begin{defn}
The $F$-translation operator is the linear operator \\ $E^{y}(\partial_{F}):
{\bf P}\rightarrow{\bf P}$ of the form:
$$E^{y}(\partial_{F})=\exp_{F}\{ y \partial_{F} \}=
\sum_{k \geq 0} \frac{y^{k} \partial_{F}^{k}}{F_{k}!},\;\;\;\;\;\;
y\in{\bf K}$$
\end{defn}
\begin{defn}
$$\forall_{p\in{\bf P}}\;\;\;\;p(x+_{F}y)=E^{y}(\partial_{F})p(x)\;\;\;\;
x,y\in{\bf K}$$
\end{defn}
 \begin{defn}
A linear operator $T:{\bf P}\rightarrow{\bf P}$ is said to be
 $\partial_{F}$-shift invariant iff
$$\forall_{y \in {\bf K}}\;\;\;\; [T,E^{y}(\partial_{F})]=
TE^{y}(\partial_{F})-E^{y}(\partial_{F})T=0$$
We shall denote by $\Sigma_{F}$ the algebra of $F$-linear
$\partial_{F}$-shift invariant operators.
\end{defn}

 \begin{defn}
Let $Q(\partial_{F})$ be a formal series in powers of $\partial_{F}$ and
$Q(\partial_{F}):{\bf P}\rightarrow{\bf P}$. $Q(\partial_{F})$ is
said to be $\partial_{F}$-delta operator iff
\renewcommand{\labelenumi}{(\alph{enumi})}
\begin{enumerate}
\item $Q(\partial_{F}) \in \Sigma_{F}$
\item $Q(\partial_{F})(x) = const \neq 0$
 \end{enumerate}
 \end{defn}

 Under quite natural specification the proofs of most statements
 might be reffered to \cite{1}(see also references therein).

 The particularities of the case considered here are revealed in
 the sequel especially in the section 4 and 5. There the scope of
 new possibilities is initiated by means of unknown before
 examples.
\begin{prop}
Let $Q(\partial_{F})$ be the $\partial_{F}$-delta operator. Then
$$\forall_{ c\in{\bf K}}\;\;\;\;Q(\partial_{F})c=0.$$
\end{prop}
\begin{prop}
Every $\partial_{F}$-delta operator reduces degree of any polynomial by one.
\end{prop}

\begin{defn}
The polynomial sequence $\{q_{n}(x)\}_{n\geq 0}$ such that \\
$deg\;q_{n}(x)=n$ and:
\renewcommand{\labelenumi}{(\arabic{enumi})}
\begin{enumerate}
\item $q_{0}(x)=1;$
\item $q_{n}(0)=0,\; n\geq 1;$
\item $Q(\partial_{F})q_{n}(x)=F_{n}q_{n-1}(x),\;\;n \geq 0$
\end{enumerate}
 is called $\partial_{F}$-basic polynomial sequence of
the $\partial_{F}$-delta operator $Q(\partial_{F})$.
\end{defn}
\begin{prop}
For every $\partial_{F}$-delta operator $Q(\partial_{F})$ there exists
the uniquely determined $\partial_{F}$-basic polynomial sequence
$\{q_{n}(x)\}_{n\geq 0}$.
\end{prop}
\begin{defn}
 A polynomial sequence $\{p_{n}(x)\}_{n\geq 0}$ ($deg\;p_{n}(x)=n$)
 is of $F$-binomial (fibonomial) type if it satisfies the condition
 $$E^{y}(\partial_{F})p_{n}(x)=p_{n}(x+_{F}y)=\sum_{k\geq 0} \binom{n}{k}_{F}
 p_{k}(x)p_{n-k}(y)\;\;\;\forall_{  y\in {\bf K}}$$
\end{defn}
\begin{thm}
The  polynomial sequence $\{p_{n}(x)\}_{n\geq 0}$ is a $\partial_{F}$-basic
 polynomial sequence of some $\partial_{F}$-delta operator $Q(\partial_{F})$
 iff it is a sequence of\\ $F$-binomial type.
\end{thm}
\begin{thm} {\em (First Expansion Theorem)}\\
Let $T \in \Sigma_{F}$ and let $Q(\partial_{F})$ be a
$\partial_{F}$-delta operator with $\partial_{F}$-basic polynomial
sequence $\{q_{n}\}_{n \geq 0}$. Then
$$T = \sum _{n \geq 0} \frac{a_{n}}{F_{n}!}Q(\partial_{F})^{n};
\quad
a_{n} = [Tq_{k}(x)]_{x=0}.$$
\end{thm}

\begin{thm}{\em (Isomorphism Theorem)}\\
Let $\Phi_{F}={\bf K}_{F}[[t]]$ be the algebra of formal {\em exp}$_{F}$
series in $t \in {\bf K}$ ,i.e.:
$$f_{F}(t)\in \Phi_{F}\;\;\;\;  iff\;\;\;\; f_{F}(t)= \sum_{k \geq 0}
\frac{a_{k}t^{k}}{F_{k}!}\;\;\; for\;\;\; a_{k}\in {\bf K},$$
and let the $Q(\partial_{F})$ be a
$\partial_{F}$-delta operator. Then $\Sigma_{F} \approx \Phi_{F}$.
 The isomorphism \\ $\phi : \Phi_{F} \rightarrow \Sigma_{F}$ is given by the
 natural correspondence:
 $$f_{F}  \left( {t} \right) = \sum_{k \geq 0}
\frac{a_{k}t^{k}}{F_{k}!}\;
\buildrel {into} \over \longrightarrow \; T_{\partial_{F}} =
\sum_{k \geq 0} \frac{a_{k}}{F_{k}!}Q(\partial_{F})^{k}. $$
\end{thm}
\begin{rem} {\em
In the algebra $\Phi_{F}$ the product is given by the fibonomial convolution,
i.e.:
$$
\left( {\;\sum_{k \geq 0} {\frac{{a_{k}} }{{F_{k}  !}}} x^{k}\;}
\right)
\left( {\;\sum_{k \geq 0} {\frac{{b_{k}} }{{F_{k}  !}}} x^{k}\;}
\right)=
\left( {\;\sum_{k \geq 0} {\frac{{c_{k}} }{{F_{k}  !}}} x^{k}\;}
\right)$$
where
$$c_{k} = \sum_{l \geq 0} \binom{k}{l}_{F}  a_{l} b_{k -l}.$$ }
\end{rem}
\begin{cor}
Operator $T \in \Sigma_{F}$ has its
inverse $T^{ - 1}\; \in \Sigma_{\psi}$   iff    $T1 \neq 0$.
\end{cor}
\begin{rem}{\em
 The $F$-translation operator
$E^{y}\left( {\partial _{F} }  \right) = \exp_{F}  \{ y\partial _{F} \} $
is invertible in $\Sigma_{F}$ but it is not a $\partial _{F} $-delta operator.
 No one of $\partial_{F}  $-delta operators $Q\left( {\partial _{F} }  \right)$
is invertible with respect to the formal series "F-product".  }
\end{rem}
\begin{cor}
Operator $R(\partial _{F} ) \in \Sigma_{F} $ is a $\partial
_{F}$-delta operator iff $a_{0} = 0$ and $a_{1} \neq 0$, where
$R(\partial _{F} )  = \sum_{n \geq 0} \frac{{a_{n}} }{F_{n}
!} Q\left( \partial _{F}  \right)^{n}$ or equivalently :
$r(0) = 0$ {\em \&} $r'(0) \neq 0$ where $r(x) =
\sum\limits_{k \geq 0}\frac{a_{k}}{F_{k}!}x^{k}\;$
is the correspondent of $R(\partial _{F})\;$ under the Iomorphism Theorem.
\end{cor}
\begin{cor}
Every $\partial_{F}  $-delta operator $Q\left( \partial
_{F}  \right)$ is a function $Q(\partial _{F})  $ according to the expansion
$$Q\left( \partial _{F}   \right) = \sum\limits_{n \geq 1}\frac{q_{n}
}{F_{n}  !} \partial _{F} ^{n}$$
This $F$-series will be called the $F$-indicator of the $Q(\partial_{F})$.
\end{cor}
\begin{rem}{\em
$\exp_{F} \{zx\}$ is the $F$-exponential generating function
for \\ $\partial _{F}  $-basic polynomial sequence $\left\{ x^{n}
\right\}_{n = 0}^{\infty}  $ of the $\partial _{F}  $ operator.}
\end{rem}

\begin{cor}
The $F $-exponential generating function for $\partial _{F}  $-basic
polynomial sequence $\left\{ p_{n} \left( {x} \right) \right\}_{n =
0}^{\infty}  $ of the $\partial _{F}  $-delta operator $Q\left( {\partial
_{F} }  \right)$ is given by the following formula
$$\sum\limits_{k \geq 0} \frac{{p_{k} \left( x \right)}}{F_{k}  !}
z^{k}\; = \exp_{F}  \{ xQ^{ - 1}\left( z \right)\}$$
where
$$Q \circ Q^{-1}=Q^{-1} \circ Q=I=id.$$
\end{cor}
 \begin{ex} \label{deltaop}

  {\em The following operators are the examples of $\partial_{F}$-delta\\
  operators:}
\renewcommand{\labelenumi}{\em (\arabic{enumi})}
\begin{enumerate} {\em
\item $\partial_{F}$;
\item $F$-difference operator $\Delta_{F}=E^{1}(\partial_{F})-I$
 such that\\
$(\Delta_{F}p)(x)=p(x+_{F}1)-p(x)$ for every $p \in {\bf P}$ ;
\item The operator $\nabla_{F}=I-E^{-1}(\partial_{F})$ defined as follows:\\
$(\nabla_{F}p)(x)=p(x)-p(x-_{F}1)$ for every $p \in {\bf P}$;
\item $F$-Abel operator: $A(\partial_{F})=\partial_{F}E^{a}(\partial_{F})=
\sum\limits_{k \geq 0}
\frac{a^{k}}{F_{k}!}\partial_{F}^{k+1}$;
\item $F$-Laguerre operator of the form: $L(\partial_{F})=\frac{\partial_{F}}
{\partial_{F}-I}=\sum\limits_{k \geq 0}\partial_{F}^{k+1}$.}
\end{enumerate}
\end{ex}
 \section{The Graves-Pincherle $F$-derivative}
\begin{defn}
The $\hat{x}_{F}$-operator is the linear map $\hat{x}_{F}:{\bf P}
\rightarrow {\bf P}$ such that\\ $\hat{x}_{F}x^{n}=\frac{n+1}{F_{n+1}
}x^{n+1} \; for \;\; n\geq 0$. ($\;[\partial_{F},\hat{x}_{F}]=id$.)
\end{defn}
\begin{defn}
A linear map {\bf '} : $\Sigma _{F} \rightarrow \Sigma _{F} $ such that   \\
$T\;${\bf '} = $T\;\hat {x}_{F}  - \hat {x}_{F}  T$ = \textbf{[}$T$, $\hat
 {x}_{F} $\textbf{]}\\
 is called the Graves-Pincherle $F$-derivative {\em \cite{9,10}}.
\end{defn}
\begin{ex} \textrm{ }
 \renewcommand{\labelenumi}{(\arabic{enumi})}
\begin{enumerate}
\item $\partial_{F}${\bf '}=$I=id$;
\item $(\partial_{F})^{n}${\bf '}=$n\partial_{F}^{n-1}$
\end{enumerate}
\end{ex}

 According to the example above the Graves-Pincherle
$F$-derivative is the formal derivative with respect to $\partial_{F}$ in
$\Sigma_{F}$ i.e., $T${\bf }'$\;(\partial_{F}) \in \Sigma_{F}$ for any $T
\in \Sigma_{F} $.
\begin{cor}\label{corfourone}
Let $t\left( {z} \right) $ be the indicator of operator $T \in \Sigma _{F} $.
Then\\ $t'\left( {z} \right)$ is the indicator of   $T${\bf '}$ \in
\Sigma _{F}$.
\end{cor}

Due to the isomorphism theorem and the Corollaries above the Leibnitz rule
 holds .
\begin{prop} \label{propfourone}
($T S$)\textbf{'} $=T$\textbf{'} $S\; + ST$\textbf{'}$\;$ ;
$T$, $S \in \; \Sigma _{F} $.
\end{prop}

As an immediate consequence of the Proposition~\ref{propfourone} we get
\begin{center}
($S  ^{n}\;$)\textbf{'}= n $S$\textbf{'}$S^{n - 1}\; \quad \forall_{ S \in
\Sigma _{F}}$.
\end{center}
From the isomorphism theorem we insert that the following is true.

\begin{prop}
$Q\left( {\partial _{F} }  \right)$ is the $\partial _{F}$-delta operator iff
 there exists invertible
$S\in \Sigma_{F}$ such that
$$Q\left( {\partial _{F} }  \right) \; = \;
\partial _{F} S.$$
\end{prop}

The Graves-Pincherle $F$-derivative notion appears very effective while
formulating expressions for $\partial _{F}  $-basic polynomial sequences
of the given $\partial _{F}  $-delta operator $Q\left( {\partial _{F}
} \right)$.
  \begin{thm} \label{lr}
{\em ($F$-Lagrange and $F$-Rodrigues formulas) \cite{1,6,12}}\\
Let $\{q_{n}\}_{n \geq 0}$ be $\partial_{F}$-basic sequence of the
delta operator $Q(\partial_{F})$, $Q(\partial_{F})=\partial_{F}P$ ($P \in
\Sigma_{F}$, invertible). Then for $n\geq0$:

\begin{enumerate}
\renewcommand{\labelenumi}{\em (\arabic{enumi})}
\item\label{one} $q_{n}(x) = Q\left( \partial _{F} \right)$\textbf{'}
 $P^{-n-1}\;x^{n}$ ;

\item\label{two} $q_{n}(x) = P^{-n} x^{n} - \frac{F_{n}}{{n}}$ ($P^{ - n}\;$)
\textbf{'}$x^{n-1};$

\item\label{three} $q_{n}(x) = \frac{{F_{n} } }{{n}}\hat {x}_{F}P^{ - n}
x^{n-1}$;

\item\label{four} $q_{n}(x) = \frac{F_{n}}{n}\hat {x}_{F}(Q\left(
{\partial _{F} }  \right)$\textbf{'} )$^{-1}
q_{n-1}(x)$  {\em ($\leftarrow$ Rodrigues $F $-formula )}.
\end{enumerate}
\end{thm}
\begin{cor}
Let $Q(\partial_{F})=\partial_{F}S$ and $R(\partial_{F})=\partial_{F}P$ be
 the $\partial_{F}$-delta operators with the $\partial_{F}$-basic sequences
 $\{q_{n}(x)\}_{n \geq 0}$ and $\{r_{n}(x)\}_{n \geq 0}$ respectively. Then:
 \renewcommand{\labelenumi}{\em (\arabic{enumi})}
\begin{enumerate}
\item $q_{n}(x)=R$\textbf{'}$(Q$\textbf{'}$)^{-1}S^{-n-1}P^{n+1}r_{n}(x),
\;\;\;n \geq 0$;
\item $q_{n}(x)=\hat{x}_{F}(PS^{-1})^{n}\hat{x}_{F}^{-1}r_{n}(x),\;\;\;
n>0$.
\end{enumerate}
\end{cor}
The formulas of the Theorem~\ref{lr} can be used to find $\partial_{F}$-basic
 sequences of the $\partial_{F}$-delta operators from the
  Example~\ref{deltaop}.
\begin{ex} \textrm{  }
\renewcommand{\labelenumi}{\em (\arabic{enumi})}
\begin{enumerate} {\em
\item The polynomials $x^{n},\;n \geq 0$ are $\partial_{F}$-basic for
 $F$-derivative $\partial_{F}$.
 \item Using Rodrigues formula in a straighford  way one can find the following first
 $\partial_{F}$-basic polynomials
 of the operator $\Delta_{F}$:\\
 $q_{0}(x)=1\\
 q_{1}(x)=x\\
 q_{2}(x)=x^{2}-x\\
 q_{3}(x)=x^{3}-4x^{2}+3x\\
 q_{4}(x)=x^{4}-9x^{3}+24x^{2}-16x\\
 q_{5}(x)=x^{5}-20x^{4}+112.5x^{3}-250x^{2}+156.5x\\
 q_{6}(x)=x^{6}-40x^{5}+480x^{4}-2160x^{3}+4324x^{2}-2605x .$
 \item Analogously to the above example we find the following first
  $\partial_{F}$-basic polynomials of the operator
 $\nabla_{F}$:\\
$q_{0}(x)=1\\
 q_{1}(x)=x\\
 q_{2}(x)=x^{2}+x\\
 q_{3}(x)=x^{3}+4x^{2}+3x\\
 q_{4}(x)=x^{4}+9x^{3}+24x^{2}+16x\\
 q_{5}(x)=x^{5}+20x^{4}+112.5x^{3}+250x^{2}+156.5x\\
 q_{6}(x)=x^{6}+40x^{5}+480x^{4}+2160x^{3}+4324x^{2}+2605x .$
 \item Using Rodrigues formula in a straighford  way one finds the following first
 $\partial_{F}$-basic polynomials of $F$-Abel operator:\\
$A^{(a)}_{0,F}(x)=1\\
A^{(a)}_{1,F}(x)=x\\
A^{(a)}_{2,F}(x)=x^{2}+ax\\
A^{(a)}_{3,F}(x)=x^{3}-4ax^{2}+2a^{2}x\\
A^{(a)}_{4,F}(x)=x^{4}-9ax^{3}+18a^{2}x^{2}-3a^{3}x .$

  \item In order to find $\partial_{F}$-basic polynomials of
   $F$-Laguerre operator $L(\partial_{F})$  we use formula (3) from Theorem~\ref{lr}:
   \begin{multline*}
   L_{n,F}(x)=\frac{F_{n}}{n}\hat{x}_{F}\left(\frac{1}{\partial_{F}-1}
   \right)^{-n}x^{n-1}=\frac{F_{n}}{n}\hat{x}_{F}(\partial_{F}-1)^{n}
   x^{n-1}=\\=\frac{F_{n}}{n}\hat{x}_{F}\sum_{k=0}^{n}(-1)^{k}\binom{n}{k}
   \partial_{F}^{n-k}x^{n-1}=\frac{F_{n}}{n}\hat{x}_{F}\sum_{k=0}^{n}
   (-1)^{k}\binom{n}{k}(n-1)^{\underline{n-k}}_{F} x^{k-1}=\\=
   \frac{F_{n}}{n}\sum_{k=1}^{n}(-1)^{k}\binom{n}{k}(n-1)^{\underline
   {n-k}}_{F} \frac{k}{F_{k}}x^{k}.
  \end{multline*}}
  \end{enumerate}
 \end{ex}

 \section{Sheffer $F$-polynomials}
\begin{defn}
A polynomial sequence $\{s_{n}\}_{n\geq 0}$ is called the sequence of Sheffer
$F$-polynomials of the $\partial_{F}$-delta operator
$Q(\partial_{F})$ iff

\renewcommand{\labelenumi}{\em (\arabic{enumi})}
\begin{enumerate}
\item $s_0(x)=const\neq 0$
\item $Q(\partial_{F})s_{n}(x)=F_{n}s_{n-1}(x);\; n\geq 0.$
\end{enumerate}
\end{defn}

\begin{prop} \label{shefprop}
Let $Q(\partial_{F})$ be $\partial_{F}$-delta operator with
$\partial_{F}$-basic polynomial sequence $\{q_{n}\}_{n  \geq 0}$. Then
$\{s_{n}\}_{n\geq 0}$ is the sequence of Sheffer $F$-polynomials of
$Q(\partial_{F})$ iff there exists an invertible $S \in \Sigma
_{F}$ such that $s_{n}(x)=S^{-1}q_{n}(x)$ for  $n\geq 0$.
We shall refer to a given labeled by  $\partial_{F}$-shift invariant
invertible operator $S$   Sheffer $F$-polynomial sequence $\{s_{n}\}
_{n\geq 0}$ as the sequence of Sheffer $F$-polynomials of the
$\partial_{F}$-delta operator $Q(\partial_{F})$ relative to $S$.
\end{prop}
\begin{thm}\label{thfourthree} {\em (Second $F $- Expansion Theorem)}\\
Let $Q\left( {\partial _{F} }  \right)$ be the $\partial _{F}  $-delta
operator $Q\left( {\partial _{F} }  \right)$ with the $\partial _{F}
$-basic polynomial sequence $\left\{ {q_{n} \left( {x} \right)}
\right\}_{n\geq 0}  $. Let $S$ be an \textit{invertible}
$\partial _{F}  $-shift invariant operator and let $\left\{ {s_{n} \left(
{x} \right)} \right\}_{n \geq 0} $ be its sequence of Sheffer $F
$-polynomials. Let $T$ be \textit{any}  $\partial
_{F}  $-shift invariant operator and let \textit{p(x)} be any polynomial.
Then the following identity holds :

\begin{center}
$\forall_{ y \in K} \wedge \; \forall_{ p \in P} \; \quad (Tp)\left( {x
+ _{F} y} \right) =\left[E^{y}(\partial_{F})p\right](x)=T \sum\limits
_{k \geq 0}\frac{s_{k} \left( y \right)}{F_{k}  !} Q\left( {\partial
_{F} }  \right)^{k}S\;  Tp\left( {x} \right)$ .
\end{center}
\end{thm}
\begin{cor}
Let ${s_{n}(x)}_{n \geq 0}$ be a sequence of Sheffer $F$-polynomials
of  a  $\partial_{F}$-delta operator $Q(\partial_{F})$ relative to $S$.Then:
$$S^{-1}=\sum_{k \geq 0}\frac{s_{k}(0)}{F_{k}!}Q(\partial_{F})^{k}.$$
\end{cor}
\begin{thm}({\em The Sheffer $F$-Binomial Theorem})\\
Let $Q(\partial_{F})$, invertible $S \in\Sigma_{F},{q_{n}(x)}
_{n \geq 0},{s_{n}(x)}_{n \geq 0}$ be as above. Then:
 \begin{center}
$$E^{y}(\partial_{F})s_{n}(x)=s_{n}(x+_{F}y)=\sum_{k \geq 0}\binom{n}{k}_{F}
s_{k}(x)q_{n-k}(y).$$
\end{center}
\end{thm}
\begin{cor}
$$s_{n}(x)=\sum_{k \geq 0}\binom{n}{k}_{F}s_{k}(0)q_{n-k}(x)$$
\end{cor}
\begin{prop}
Let $Q\left( {\partial _{F} }  \right)$ be a $\partial _{F}  $-delta
operator. Let $S$ be an invertible $\partial _{F
} $-shift invariant operator. Let $\left\{ {s_{n} \left( {x} \right)}
\right\}_{n \geq 0} $ be a polynomial sequence. Let
\begin{center}
$\forall_{ a \in K }\wedge \; \forall_{ p \in P} \quad E^{a}\left(
{\partial _{F} }  \right)p\left( {x} \right) = \sum\limits_{k \ge 0}
{\frac{{s_{k} \left( {a} \right)}}{{F_{k}  !}}} Q\left( {\partial
_{F} }  \right)^{k}S_{\partial _{F} }  \;p\left( {x} \right)$ .
\end{center}
Then the polynomial sequence $\left\{ {s_{n} \left( {x} \right)} \right\}_{n
\geq 0}$ is the sequence of Sheffer  $F $-polynomials of the
$\partial _{F}  $-delta operator $Q\left( {\partial _{F} }  \right)$
relative to $S$.
\end{prop}
\begin{prop}
Let $Q\left( {\partial _{F} }  \right)$and  $S$ be as above. Let
\textit{q(t)} and s\textit{(t)}  be the
indicators of $Q\left( {\partial _{F} }  \right)$ and $S$ operators.
Let \textit{q$^{-1}$(t )}
be the inverse $F $-exponential formal power series inverse to
\textit{q(t)}. Then the $F $-exponential generating function of Sheffer
$F $-polynomials sequence $\left\{ {s_{n} \left( {x} \right)} \right\}_{n
\geq 0}$ of $Q\left( {\partial _{F} }  \right)$ relative
to $S\;$is given by
\[
\;\sum\limits_{k \ge 0} {\frac{{s_{k} \left( {x} \right)}}{{F_{k}
!}}} z^{k}\; = \;\left(s\left( {q^{ - 1}\left( {z} \right)}
\right)\right)^{-1}\;\exp_{F}  \{ xq^{ - 1}\left( {z} \right)\}.
\]
\end{prop}
\begin{prop}
A sequence $\left\{ {s_{n} \left( {x} \right)} \right\}_{n \geq 0}$
is the sequence of Sheffer  $F $-polynomials of the $\partial _{F}
$-delta operator $Q\left( {\partial _{F} }  \right)$ with the $\partial
_{F}  $-basic polynomial sequence $\left\{ {q_{n} \left( {x} \right)}
\right\}_{n \geq 0}$ iff
\[
s_{n} \left( {x + _{F}  y} \right)=
\sum\limits_{k \ge 0} \binom{n}{k}_{F}  s_{k} \left( {x} \right)q_{n - k} \left(
{y} \right).
\]
for all $y \in {\bf K}$
\end{prop}

\begin{ex}{\em
Hermite $F$-polynomials are Sheffer $F$-polynomials of the\\ $\partial_{F}$
-delta operator $\partial_{F}$ relative to invertible $S \in \Sigma_{F}$ of
the form \\ $S=\exp_{F} \{ \frac{a \partial_{F}^{2}}{2}\}$. One can get
them by formula (see Proposition~\ref{shefprop} ):
$$H_{n,F}(x)=S^{-1}x^{n}=\sum\limits_{k \geq 0}\frac{(-a)^{k}}{2^{k}F_{k}!}
n^{\underline{2k}}_{F}x^{n-2k}.$$}
\end{ex}
\begin{ex}   {\em
 Let $S=(1-\partial_{F})^{-\alpha- 1}$. The Sheffer $F$-polynomials of\\
 $\partial_{F}$-delta operator $L(\partial_{F})=\frac{\partial_{F}}
 {\partial_{F}-1}$  relative to $S$ are Laguerre $F$-polynomials of order
 $\alpha$ . By Proposition~\ref{shefprop} we have
 $$L^{(\alpha)}_{n,F}=(1-\partial_{F})^{\alpha+1}L_{n,F}(x),$$
From the above formula and using Graves-Pincherle $F$-derivative
we get
$$L^{(\alpha)}_{n,F}(x)=\sum\limits_{k \geq 0}\frac{F_{n}!}{F_{k}!}\binom
{\alpha+n}{n-k}(-x)^{k}$$ for $\alpha \neq -1$}.
\end{ex}
\begin{ex} {\em
Bernoullie's $F$-polynomials of order 1 are Sheffer $F$-polynomials of \\$\partial_{F}$
-delta operator $\partial_{F}$ related to invertible $S=\left(\frac
{\exp_{F}\{\partial_{F}\}-I}{\partial_{F}}\right)^{-1}$. Using\\
Proposition~\ref{shefprop} one arrives at
 \begin{multline*}
 B_{n,F}(x)=S^{-1}x^{n}=\sum_{k \geq 1}\frac{1}{F_{k}!}\partial_{F}
^{k-1}x^{n}=\sum_{k \geq 1}\frac{1}{F_{k}}\binom{n}{k-1}_{F}x^{n-k+1}=\\=
\sum_{k \geq 0}\frac{1}{F_{k+1}}\binom{n}{k}_{F}x^{n-k}
\end{multline*}}
\end{ex}

\begin{thm}\label{recrel}
{\em (Reccurence relation for Sheffer $F$-polynomials)}\\
Let $Q,S,\{s_{n}\}_{n\geq 0} $ be as above. Then the following reccurence
formula holds:
$$
s_{n+1}(x)=\frac{F_{n+1}}{n+1}\left[\hat{x}_{F}-\frac{S'}{S}\right]
\left[Q(\partial_{F})'\right]^{-1}s_{n}(x);\; n\geq 0.
$$
\end{thm}
\begin{ex}
{\em The reccurence formula for the Hermite $F$-polynomials is:}
$$H_{n+1,F}(x)=\hat{x}_{F}H_{n,F}(x)-\hat{a}_{F}F_{n}H_{n-1,F}(x)$$
\end{ex}
\begin{ex}
{\em The reccurence relation for the Laguerre $F$-polynomials is:}
\begin{multline*}
L_{n+1,F}^{(\alpha)}(x)=-\frac{F_{n+1}}{n+1}[\hat{x}_{F}-(\alpha+1)(1-\partial_{F})^{-1}
](\partial_{F}-1)^{2}L_{n,F}^{(\alpha)}(x)\\=\frac{F_{n+1}}{n+1}[\hat{x}_{F}(\partial_{F}-1)
+\alpha+1]L_{n,F}^{(\alpha+1)}(x).
\end{multline*}
\end{ex}
  \section{The Spectral Theorem}

We shall now define a natural inner product associated with the sequence
$\{s_{n}\}_{n\geq 0}$  of Sheffer $F$-polynomials of the
$\partial_{F}$-delta operator $Q(\partial_{F})$ relative to $S$.

\begin{defn}
Let $Q,S,\{s_{n}\}_{n\geq 0} $ be as above. Let $W$ be umbral operator:\\
$W:s_{n}(x)\rightarrow x^{n}$ ( and linearly extented). We define the
following bilinear form:
$$(f(x),g(x))_{F}:=[(Wf)(Q(\partial_{F}))Sg(x)]_{x=0};\;\;
f,g\in{\bf P}.
$$
\end{defn}

\begin{prop} \label{propthree} { \em \cite{6}}
The bilinear form over reals defined above is a positive definite
 inner product such
that:
$$(s_{n}(x),s_{k}(x))_{F}=F_{n}!\delta_{n,k}.$$
We shall call this scalar praduct the natural inner product associated with
the sequence $\{s_{n}\}_{n\geq 0}$ of Sheffer $F$-polynomials. Unitary
space $({\bf P},(\textrm{ },\textrm{ })_{F})$ can be completed to the
unique Hilbert space ${\bf H}=\overline{{\bf P}}$.
\end{prop}

\begin{thm} {\em (Spectral Theorem)}\\
Let $\{s_{n}\}_{n\geq 0}$ be the sequence of Sheffer $F$-polynomials
relative to the \\ $\partial_{F}$-shift invariant invertible operator $S$
for the $\partial_{F}$-delta operator $Q(\partial_{F})$ with
$\partial_{F}$-basic polynomial sequence $\{q_{n}\}_{n\geq 0}$. Then
there exists a unique operator $A_{F}:{\bf H}\rightarrow{\bf H}$ of the
form
$$A_{F}=\sum \limits_{k\geq 1} \frac{u_{k}+\hat{v}_{k}(x)}{F
_{k-1}!}Q(\partial_{F})^{k}$$  with the following properties:
\renewcommand{\labelenumi}{(\alph{enumi})}
\begin{enumerate}
\item $A$ is self adjoint;
\item The spectrum of $A$ consists of $n\in{\bf N}$ and $As_{n}=ns_{n}$ for
$n\geq 0$;
\item Quantities $u_{k}$ and $\hat{v}_{k}(x)$ are calculated according to
$$u_{k}=-[(\log{S})' \hat{x}_{F}^{-1}q_{k}(x)]_{x=0} \quad \quad
\hat{v}_{F}(x)=\hat{x}_{F}\left[\frac{d}{dx}q_{k}(x)\right]_{x=0}$$
{\em Proof: see \cite{1}.}
\end{enumerate}
\end{thm}
 \section{The first elementary examples of $F$-polynomials}
  \renewcommand{\labelenumi}{(\arabic{enumi})}
\begin{enumerate}
\item Here are the examples of Laguerre $F$-polynomials of order $\alpha=-1$:
\\ \textrm{}\\
$L_{0,F}(x)=1\\ \textrm{ } \\
L_{1,F}(x)=-x\\ \textrm{ } \\
L_{2,F}(x)=x^{2}-x\\ \textrm{ } \\
L_{3,F}(x)=-x^{3}+4x^{2}-2x\\ \textrm{ } \\
L_{4,F}(x)=x^{4}-9x^{3}+18x^{2}-6x\\  \textrm{ } \\
L_{5,F}(x)=-x^{5}+20x^{4}-905x^{3}+1280x^{2}-30x\\  \textrm{ } \\
L_{6,F}(x)=x^{6}-40x^{5}+400x^{4}-1200x^{3}+1200x^{2}-
240x\\  \textrm{ } \\
L_{7,F}(x)=-x^{7}+78x^{6}-1560x^{5}+
10400x^{4}-23400x^{3}+18720x^{2}-\\
 \textrm{}\\\;\;\;\;\;\;-3120x\\ \textrm{}\\
L_{8,F}(x)=x^{8}-147x^{7}+5733x^{6}-76440x^{5}+
382200x^{4}-687960x^{3}+\\ \;\;\;\;+458640x^{2}-65520x$
\item Here are the examples of Laguerre $F$-polynomials of order $\alpha=1$:
\\ \textrm{}\\
$L^{(1)}_{0,F}(x)=1\\ \textrm{ } \\
L^{(1)}_{1,F}(x)=-x+2\\ \textrm{ } \\
L^{(1)}_{2,F}(x)=x^{2}-3x+3\\ \textrm{ } \\
L^{(1)}_{3,F}(x)=-x^{3}+8x^{2}-12x+8\\ \textrm{ } \\
L^{(1)}_{4,F}(x)=x^{4}-15x^{3}+60x^{2}-60x+30\\  \textrm{ } \\
L^{(1)}_{5,F}(x)=-x^{5}+30x^{4}-225x^{3}+600x^{2}-450x+240\\
 \textrm{ } \\
L^{(1)}_{6,F}(x)=x^{6}-56x^{5}+840x^{4}-4200x^{3}+8400x^{2}-
5040x+1680$
\item  Here we give some examples of the Bernoullie's $F$-polynomials of
order 1:\\
\textrm{}\\
$B_{0,F}(x)=1\\ \textrm{ } \\
B_{1,F}(x)=x+1\\ \textrm{ } \\
B_{2,F}(x)=x^{2}+x+\frac{1}{2}\\ \textrm{ } \\
B_{3,F}(x)=x^{3}+2x^{2}+x+\frac{1}{3}\\ \textrm{ } \\
B_{4,F}(x)=x^{4}+3x^{3}+3x^{2}+x+\frac{1}{5}\\  \textrm{ } \\
B_{5,F}(x)=x^{5}+5x^{4}+\frac{15}{2}x^{3}+5x^{2}+x+\frac{1}{8}\\ \textrm{ } \\
B_{6,F}(x)=x^{6}+8x^{5}+20x^{4}+20x^{3}+8x^{2}+x+\frac{1}{13}\\  \textrm{ } \\
B_{7,F}(x)=x^{7}+13x^{6}+52x^{5}+\frac{260}{3}x^{4}+52x^{3}+13x^{2}+x+
\frac{1}{21}\\ \textrm{}\\
B_{8,F}(x)=x^{8}+21x^{7}+\frac{273}{2}x^{6}+364x^{5}+364x^{4}+\frac{273}{2}
x^{3}+21x^{2}+x+\frac{1}{36}\\ \textrm{}\\
B_{9,F}(x)=x^{9}+34x^{8}+357x^{7}+1547x^{6}+\frac{12376}{5}x^{5}+1547x^{4}+
357x^{3}+\\ \textrm{}\\+34x^{2}+x+\frac{1}{55}$
\end{enumerate}
\begin{rem}{ \em
Let us observe that analogously to the ordinary case
$F$-polynomials ,such as Abel, Laguerre or Bernoullie's
$F$-polynomials may have
coefficients which are integer numbers ($F$-Abel, $F$-Laguerre)
and non-integer rationals ($F$-Bernoulli).\\
To see that recall for example the formula for Laguerre
$F$-polynomials of order -1
\\($F$-basic):
$$L_{n,F}(x)=\frac{F_{n}}{n}\sum_{k=1}^{n}(-1)^{k}\binom{n}{k}(n-1)^{\underline
   {n-k}}_{F} \frac{k}{F_{k}}x^{k}$$
and the one for $F$-Laguerre of order $\alpha \neq -1$ ($F$-Sheffer):
$$L^{(\alpha)}_{n,F}(x)=\sum\limits_{k \geq 0}\frac{F_{n}!}{F_{k}!}\binom
{\alpha+n}{n-k}(-x)^{k}.$$ Because Fibonomial coefficients are
integers the second formula gives
 us polynomials with integer coefficients. It is easy to verify that
 $F$-basic \\Laguerre polynomials do have this property too.\\
 Finally let $p \in  {\bf P}$ while $a_{k}$ denote coefficient of this
 polynomial $p$ at $x^{k}$,i.e.
 $$p(x)=\sum_{k \geq 0}a_{k}x^{k}.$$
 Consider now the Bernoullie's $F$-polynomials of order 1.
 Because of the symmetry  of $\binom{n}{k}_{F}$ and some known
 divisibility properties
 of Fibonacci numbers \cite{5,8}  for Bernoullie's $F$-polynomial $B_{n,F}(x)$ we have
 $$a_{n-k}=a_{k+1}$$ for $k=0,1,..,\left[\frac{n}{2}\right]$. Moreover from
 formula for these polynomials it comes that $$a_{0}=\frac{1}{F_{n+1}}.$$ \\
 Observe now that  coefficients of Abel $F$-polynomials are integer numbers,
  so we may expect  now that these polynomials enumerate some combinatorial
  objects like those of the now classical theory of binomial enumeration
   (see \cite{11}).

 }
\end{rem}
{\bf {Acknowledgements}}

I would like to thank to Prof. A.K.Kwa\'sniewski for his remarks and
 guideness.

 AMS Classification numbers: 11C08, 11B37, 47B47

\end{document}